\documentclass[12pt]{article}
\usepackage{amsmath}
\usepackage{amssymb}
\usepackage[cp1251]{inputenc}
\usepackage[russian]{babel}
\newcommand{\il}[2]{\int\limits_{#1}^{#2}}
\newcommand{\ilp}[1]{\int\limits_{#1}^{+\infty}}

\newcommand{\ph}{\phantom{a}}
\newcommand{\phh}{\phantom{aaa}}

\newcommand{\sist}[2]{\left\{
\begin{array}{l}
{#1}\\
\ph\\
{#2}
\end{array}
\right.}

\newcommand{\liml}[1]{\lim\limits_{#1}}

\begin{document}

MSC 34D20

\vskip 20pt

\centerline{\bf Global solvability, stability and oscillation  criteria for systems of two}
\centerline{\bf first-order pseudo linear ordinary differential equations}

\vskip 10 pt

\centerline{\bf G. A. Grigorian}

\vskip 10 pt

\centerline{0019 Armenia c. Yerevan, str. M. Bagramian 24/5}
\centerline{Institute of Mathematics of NAS of Armenia}
\centerline{E - mail: mathphys2@instmath.sci.am, \ph phone: 098 62 03 05, \ph 010 35 48 61}

\vskip 20 pt

\noindent
Abstract. In this paper we use the Riccati equation method with other ones to establish   global solvability, stability and oscillation criteria for a class of two dimensional  nonlinear systems of ordinary differential equations, which is a generalization of wide classes of second order nonlinear ordinary differential equations, studied by many authors. The applicability of some of these criteria are illustrated by examples.

\vskip 20 pt

\noindent
Key words: the system of pseudo linear equations, the Riccati equation, global solvability, the van der Pol's equation, the Duffing's equation, the equation of Emden-Fowler.

\vskip 20 pt
{\bf  1. Introduction}. Let $P(t,u,v), \ph Q(t,u,v), \ph R(t,u,v), \ph S(t,u,v), \ph F(t,u,v)$ and \linebreak $G(t,u,v)$ be real-valued locally integrable  by $t$ and continuous by $u$ and $v$  functions on $[t_0,+\infty) \times (\mathbb{R}\backslash \{0\})\times  (\mathbb{R}\backslash \{0\})$, where $\mathbb{R}\equiv (-\infty,+\infty)$. Consider the system of two first-order pseudo linear (i.e. like of linear) ordinary differential equations
$$
\sist{\phi'=P(t,\phi,\psi)\phi + Q(t,\phi,\psi)\psi + F(t,\phi,\psi),}{\psi' = R(t,\phi,\psi)\phi + S(t,\phi,\psi)\psi + G(t,\phi,\psi), \ph t \ge t_0.} \eqno (1.1)
$$
In the particular case, when $F(t,u,v) = G(t,u,v) \equiv 0$ we have a "pseudo homogeneous" \linebreak case of (1.1):
$$
\sist{\phi'=P(t,\phi,\psi)\phi + Q(t,\phi,\psi)\psi,}{\psi' = R(t,\phi,\psi)\phi + S(t,\phi,\psi)\psi, \ph t \ge t_0.} \eqno (1.2)
$$
By a solution of the system (1.1) (of the system (1.2)) on $[t_1,t_2) \ph (t_0 \le t_1 < t_2 \le +\infty)$ we mean an ordered pair $(\phi(t), \psi(t))$ of  absolutely continuous  on $[t_1,t_2)$ functions $\phi(t)$ and $\psi(t)$, satisfying (1.1) ((1.2)) almost everywhere on $[t_1,t_2)$. According to the general theory of normal systems of ordinary differential equations for every $\alpha, \beta \in \mathbb{R}\backslash \{0\}$ and $t_1 \ge t_0$ there exists $t_2 > t_1 \ph (t_2 \le + \infty)$ and a solution $(\phi(t), \psi(t))$ of the system (1.1) ((1.2)) on $[t_1,t_2)$ with $\phi(t_1) = \alpha, \ph \psi(t_1) = \beta.$   A great interest from the point of view of qualitative theory of differential equations represents the case $t_2 = +\infty$. In this case we say that the system (1.1) ((1.2)) is global solvable.
Indicate some well (and not so well) known non linear second order ordinary differential equations, which are reducible to the systems of type (1.1) or (1.2):

\noindent
1. the van der Pol equation with parametric excitation (see [1], p. 333):
$$
\phi'' + \varepsilon(\phi^2 - 1)\phi' + (1 + \beta \cos t) \phi = 0; \eqno (1.3)
$$

\noindent
2. the van der Pol equation  (see [1], p. 234):
$$
\phi'' + \varepsilon(\phi^2 - 1)\phi' + \phi = \Gamma\cos \omega t ; \eqno (1.4)
$$

\noindent
3. the equation of oscillation of an electronic contour with parametric excitation having an electronic lamp (see [2], p. 221):
$$
\phi'' +2(\lambda_1 + \lambda_2\phi^2)\phi' +\omega^2(1 - h \cos\nu t)\phi = 0; \eqno (1.5)
$$

\noindent
4. the van der Pol-Matheu equation for the dynamics of dust grain charges in dusty plasma (see [3]):
$$
\phi'' - (\alpha - \beta\phi^2)\phi' +\omega_0^2(1 + h \cos \gamma t) \phi = 0; \eqno (1.6)
$$

\noindent
5. the Lienard's equation (see [1], p. 387):
$$
\phi'' + f(t,\phi')\phi' + g(\phi) = 0; \eqno (1.7)
$$

\noindent
6. the equation of a pendulum with a light (see [1], p. 331):
$$
\phi'' + \biggl(-\frac{g}{2} + \frac{\varepsilon}{a} \cos \omega t\biggr)\sin \phi  = 0; \eqno (1.8)
$$

\noindent
7. the equation of a pendulum (see [2], p. 34):
$$
(m l^2(t) \phi')' + g l(t) \sin \phi = 0; \eqno (1.9)
$$

\noindent
8. the equation of a pendulum with bob of mass $m$ and rigid suspension of length  $a$, hanging from a support, which is constrained to move  with vertical and horizontal displacement $\zeta(t)$ and $\eta(t)$ respectively (see [1], p. 334);
$$
a \phi'' + (g + \zeta''(t))\sin\phi + \eta''(t) \cos \phi = 0; \eqno (1.10)
$$

\noindent
9. the Duffing's equation (see [1], p 223):
$$
\phi'' + k\phi' + \alpha\phi + \beta\phi^3 = \Gamma\cos \omega t; \eqno (1.11)
$$

\noindent
10. the equation of the damped pendulum with periodic forcing of the pivot (see [1], p. 504):
$$
\phi'' + \sin \phi = \varepsilon(\gamma \sin t \sin \phi - k \phi'); \eqno (1.12)
$$

\noindent
11. the Reyligh's equation (see [1], p. 197):
$$
\phi'' + \varepsilon \Bigl(\frac{1}{3}(\phi')^3 - \phi'\Bigr) + \phi = 0; \eqno (1.13)
$$

\noindent
12. the equation for the relativistic perturbation of a planetary orbit (see [1], p. 218):
$$
\phi'' + \phi = k(1 + \varepsilon \phi^2); \eqno (1.14)
$$

\noindent
13. the nonlinear equation with the $p$-Laplasian (see [4])
$$
(\Phi_p(x'))' + \frac{2(p-1)}{t} \Phi(x') + a(t) g(x) =0, \ph t \ge t_0,
$$
$p > 1, \ph \Phi_p(y) = |y|^{p-2} y, \ph a(t) > 0, \ph x g(x) > 0, \ph x \ne 0;$

\noindent
14. (see [6])
$$
(|y'|^{\alpha -1} y')' + q(t) |y|^{\beta -1} y = 0;
$$

\noindent
15. (see [7])
$$
(a(t) \psi(x) K(x'))' + p(t) K(x') + q(t) f(x) g(x') = 0,
$$
$a(t) > 0, \ph x f(x) > 0, \ph x \ne 0, \ph K^2(y) < \gamma_1 y K(y), \ph q(t) > 0, \ph 0 < c \le g(x'), \ph \frac{f(x)}{x} > \gamma_2 > 0, \ph x \ne 0.$

In this paper we prove  global solvability, stability and oscillation  criteria for the systems (1.1) and (1.2). We illustrate the applicability of some of these criteria by examples.

{\bf 2. Auxiliary propositions}. Let $(\phi(t), \psi(t))$ be a solution of the system (1.1) on $[t_0, t_1) \ph (t_1 \le +\infty).$ We can interpret $\phi(t)$ as a solution of the linear equation
$$
\phi' = W(t)\phi + U(t), \phh t\in [t_0,t_1)
$$
where $W(t) \equiv P(t,\phi(t),\psi(t)), \ph U(t) = G(t,\phi(t),\psi(t))\psi(t) + F(t,\phi(t),\psi(t)) \psi(t)), \linebreak t\in~ [t_0,t_1)$. Then by Cauchi formula we have
$$
\phi(t) = \exp\biggl\{\il{t_0}{t} P(\tau,\phi(\tau),\psi(\tau))d\tau\biggr\}\biggl[\phi(t_0) + \il{t_0}{t}\exp\bigg\{-\il{t_0}{\tau}P(s;\phi(s);\psi(s))d s\biggr\}\times \phantom{aaaaaaaaaa}
$$
$$
\phantom{aaaaaaaaaaaa}\times\biggl(Q(\tau,\phi(\tau),\psi(\tau))\psi(\tau) + F(\tau,\phi(\tau),\psi(\tau))\biggr) d\tau\biggr], \phh t\in [t_0,t_1). \eqno (2.1)
$$
Analogously for $\psi(t)$ we obtain
$$
\psi(t) = \exp\biggl\{\il{t_0}{t} S(\tau,\phi(\tau),\psi(\tau))d\tau\biggr\}\biggl[\psi(t_0) + \il{t_0}{t}\exp\bigg\{-\il{t_0}{\tau}S(s,\phi(s),\psi(s))d s\biggr\}\times \phantom{aaaaaaaaaa}
$$
$$
\phantom{aaaaaaaaaaaa}\times\biggl(R(\tau,\phi(\tau),\psi(\tau))\phi(\tau) + G(\tau,\phi(\tau),\psi(\tau))\biggr) d\tau\biggr], \phh t\in [t_0,t_1). \eqno (2.2)
$$
Substitute the right hand part of the last equality into (2.1). After some simplifications we obtain
$$
\phi(t) = v_1(t) + \il{t_0}{t} K_1(t,\zeta)\phi(\zeta) d \zeta, \phh t\in [t_0,t_1), \eqno (2.3)
$$
where
$v_1(t)\equiv \phi(t_0)\exp\biggl\{\il{t_0}{t}P(\tau,\phi(\tau),\psi(\tau))d\tau\biggr\} +\\ \phantom{aaaaaaaaaaaaaaaaaaaaaa}+\psi(t_0)\il{t_0}{t}\exp\biggl\{\il{\tau}{t}P(s,\phi(s),\psi(s))d s + \il{t_0}{\tau}S(s,\phi(s),\psi(s))d s\biggr\}+$
$$
\il{t_0}{t}\exp\biggl\{\il{\tau}{t}P(s,\phi(s),\psi(s))d s\biggr\}\biggl[F(\tau,\phi(\tau),\psi(\tau)) + Q(\tau,\phi(t),\psi(\tau))\il{t_0}{\tau}G(s,\phi(s),\psi(s))d s\biggr] d \tau,
$$
$$
K_1(t,\zeta)\equiv R(\zeta,\phi(\zeta),\psi(\zeta))\il{\zeta}{t}\exp\biggl\{\il{\tau}{t}P(s,\phi(s),\psi(s))ds + \il{\zeta}{\tau}S(s,\phi(s),\psi(s))d s\biggr\}\times \phantom{aaaaaaaaaaaaaaaaaa}
$$
$$
\phantom{aaaaaaaaaaaaaaaaaaaaaaaaaaaaaaaaa} \times Q(\tau,\phi(\tau),\psi(\tau)) d \tau, \phh t, \zeta \in [t_0,t), \phh \zeta \le t.
$$
By similar way  from (2.1) and (2.2)   we can obtain  the equality
$$
\psi(t) = v_2(t) + \il{t_0}{t} K_2(t,\zeta)\psi(\zeta) d \zeta, \phh t\in [t_0,t_1), \eqno (2.4)
$$
where
$v_2(t)\equiv \psi(t_0)\exp\biggl\{\il{t_0}{t}S(\tau,\phi(\tau),\psi(\tau))d\tau\biggr\} +\\ \phantom{aaaaaaaaaaaaaaaaaaaaaa}+\phi(t_0)\il{t_0}{t}\exp\biggl\{\il{\tau}{t}S(s,\phi(s),\psi(s))d s + \il{t_0}{\tau}P(s,\phi(s),\psi(s))d s\biggr\}+$
$$
\il{t_0}{t}\exp\biggl\{\il{\tau}{t}S(s,\phi(s),\psi(s))d s\biggr\}\biggl[G(\tau,\phi(\tau),\psi(\tau)) + R(\tau,\phi(\tau),\psi(\tau))\il{t_0}{\tau}F(s,\phi(s),\psi(s))d s\biggr] d \tau,
$$
$$
K_2(t,\zeta)\equiv Q(\zeta,\phi(\zeta),\psi(\zeta))\il{\zeta}{t}\exp\biggl\{\il{\tau}{t}S(s,\phi(s),\psi(s))ds + \il{\zeta}{\tau}P(s,\phi(s),\psi(s))d s\biggr\}\times \phantom{aaaaaaaaaaaaaaaaaa}
$$
$$
\phantom{aaaaaaaaaaaaaaaaaaaaaaaaaaaaaaaaa} \times R(\tau,\phi(\tau),\psi(\tau)) d \tau, \phh t, \zeta \in [t_0,t), \phh \zeta \le t.
$$

Let $f_k(t), \ph g_k(t), \ph h_k(t), \ph k=1,2$ be real-valued locally integrable functions on $[t_0,+\infty)$. Consider the Riccati equations
$$
y' + f_k(t) y^2 + g_k(t) y + h_k(t) = 0, \phh t\ge t_0. \eqno (2.5_k)
$$
$k=1,2$ and the differential inequalities
$$
\eta + f_k(t) \eta^2 + g_k(t) \eta + h_k(t) \ge 0, \phh t\ge t_0. \eqno (2.6_k)
$$
$k=1,2$.
By a solution of Eq. $(2.5_k)$ (of inequality $(2.6_k)$) on $[t_1,t_2) \ph (t_0 \le t_1 < t_2 \le +\infty$) we mean an absolutely continuous function $y(t) \ph (\eta(t))$ on $[t_1,t_2)$, satisfying $(2.5_k) \ph ((2.6_k))$ almost everywhere on $[t_1,t_2), \ph k =1,2.$.

{\bf Remark 2.1.} {\it Every solution  of Eq. $(2.5_2)$ on $[t_0,t_1)$ is also a solution of the inequality $(2.6_2)$ on $[t_0,t_1)$.}

{\bf Remark 2.2.} {\it If $f_1(t) \ge 0, \ph t\in [t_0,t_1)$, then every solution of the linear equation
$$
\zeta' + g_1(t)\zeta + h_1(t) = 0, \phh t\in [t_0,t_1)
$$
is also a solution of the inequality $(2.6_1)$ on $[t_0,t_1)$.}

{\bf Theorem 2.1.} {\it  Let $y_2(t)$ be a solution of Eq. $(2.5_2)$ on $[t_0,\tau_0) \ph (t_0 < \tau_0 \le +\infty)$ and let $\eta_1(t)$ and $\eta_2(t)$ be solutions of the inequalities $(2.6_1)$ and $(2.6_2)$  respectively on $[t_0,\tau_0)$ such that $y_2(t_0) \le \eta_k(t_0) \ph k=1,2.$ In addition let the following conditions be satisfied: $f_1(t) \ge 0, \ph \gamma - y_2(t_0) + \il{t_0}{t}\exp\biggl\{\il{t_0}{\tau}[f_1(s)(\eta_1(s) + \eta_2(s)) + g_1(s)]ds\biggr\}\biggl[(f_2(\tau) - f_1(\tau))^2 y_2^2(\tau) + (g_2(\tau) - g_1(\tau)) y_2(\tau) + h_2(\tau) - h_1(\tau)\biggr] d \tau \ge 0, \ph t\in [t_0,\tau_0)$ for some $\gamma \in [y_2(t_0), \eta_1(t_0)]$. Then Eq. $(2.5_1)$ has a solution $y_1(t)$ on $[t_0,\tau_0)$ with $y_1(t_0) \ge \gamma$ and $y_1(t) \ge y_2(t), \ph t\in [t_0,\tau_0)$}

Proof. By analogy of the proof of Theorem 3.1 from [7].

In the system (1.2) substitute
$$
\psi = y \phi
$$
We obtain
$$
\sist{\phi' = [P(t;\phi;\psi) + y Q(t,\phi,\psi)]\phi + F(t,\phi,\psi),}{[y'\phi + y(t)\{[B(t,\phi,\psi) + y(t)Q(t,\phi,\psi)]\psi + F(t,\phi,\psi)\} = R(t,\phi,\psi) \phi + G(t,\phi,\psi), \ph t\ge t_0,}
$$
where $B(t,\phi,\psi)\equiv P(t,\phi,\psi) - S(t,\phi,\psi)$. It follows from here that if $(\phi(t), \psi(t))$ is a solution of the system (1.1) on $[t_0,t_1)$ with $\phi(t) \ne 0, \ph t\in [t_0,t_1)$ and $y(t)$ is the solution of the Riccati equation
$$
y' + Q(t,\phi(t),\psi(t)) y^2 + \Biggl(B(t,\phi(t),\psi(t)) + \frac{F(t,\phi(t),\psi(t))}{\phi(t)}\Biggr) y - \phantom{aaaaaaaaaaaaaaaaaaaaaaaaaaa}
$$
$$
\phantom{aaaaaaaaaaaaaaaaaaaaa}-R(t,\phi(t),\psi(t)) - \frac{G(t,\phi(t),\psi(t))}{\phi(t)} = 0 \eqno (2.7)
$$
on $[t_0,t_1)$ with $y(t_0) = \frac{\psi(t_0)}{\phi(t_0)}$ then
$$
{\left\{
\begin{array}{l}
\phi(t) = \phi(t_0)\exp\biggl\{\il{t_0}{t}\bigl[P(\tau,\phi(\tau),\psi(\tau)) + \frac{F(\tau,\phi(\tau),\psi(\tau))}{\phi(\tau)}+\\
\phantom{aaaaaaaaaaaaaaaaaaaaaaaaaaaaaaaaaaa} + y(\tau)Q(\tau;\phi(\tau);\psi(\tau))\bigr]d\tau\biggr\}, \ph \mbox{or}\\
\ph\\
\phi(t) =\exp\biggl\{\il{t_0}{t}\bigl[P(\tau;\phi(\tau),\psi(\tau)) +  y(\tau)Q(\tau,\phi(\tau),\psi(\tau))\bigr]d\tau\biggr\}\biggl[\phi(t_0) + \phantom{aaaaaaaa}\\
+ \il{t_0}{t}\exp\biggl\{- \il{t_0}{\tau}\bigl[P(s,\phi(s),\psi(s)) + y(s)Q(s,\phi(s),\psi(s))\bigr]ds \biggr\}F(\tau,\phi(\tau),\psi(\tau)) d \tau\biggr] \\
\ph\\
\psi(t) = y(t)\phi(t), \ph t\in [t_0,t_1)
\end{array}
\right.} \eqno(2.8)
$$


Analogously the substitution
$$
\phi = z\psi
$$
in (1.1) implies that if $(\phi(t),\psi(t))$ is a solution of the system (1.1) on $[t_0,t_1)$ with $\psi(t) \ne~ 0, \ph t\in [t_0,t_1)$ and $z(t)$ is the solution of the Riccati equation
$$
y' + R(t,\phi(t),\psi(t)) y^2 - \Biggl(B(t,\phi(t),\psi(t)) + \frac{G(t,\phi(t),\psi(t))}{\phi(t)}\Biggr) y - \phantom{aaaaaaaaaaaaaaaaaaaaaaaaaaa}
$$
$$
\phantom{aaaaaaaaaaaaaaaaaaaaa}-Q(t,\phi(t),\psi(t)) - \frac{F(t,\phi(t),\psi(t))}{\phi(t)} = 0 \eqno (2.9)
$$

on $[t_0,t_1)$ with  $z(t_0) = \frac{\phi(t_0)}{\psi(t_0)}$ then

$$
{\left\{
\begin{array}{l}
\psi(t) = \psi(t_0)\exp\biggl\{\il{t_0}{t}\bigl[S(\tau,\phi(\tau),\psi(\tau)) + \frac{G(\tau,\phi(\tau),\psi(\tau))}{\psi(\tau)}+\\
\phantom{aaaaaaaaaaaaaaaaaaaaaaaaaaaaaaa} + z(\tau)R(\tau,\phi(\tau),\psi(\tau))\bigr]d\tau\biggr\}, \ph \mbox{or}\\
\ph\\
\psi(t)=\exp\biggl\{\il{t_0}{t}\bigl[S(\tau;\phi(\tau);\psi(\tau)) + z(\tau)QR(\tau,\phi(\tau),\psi(\tau))\bigr]d\tau\biggr\}\biggl[\psi(t_0) + \phantom{aaaa} \\
\il{t_0}{t}\exp\biggl\{- \il{t_0}{\tau}\bigl[S(s,\phi(s),\psi(s)) + z(s)R(s,\phi(s),\psi(s))\bigr]ds \biggr\}G(\tau,\phi(\tau),\psi(\tau)) d \tau\biggr] \\
\ph\\
\psi(t) = y(t)\phi(t), \ph t\in [t_0,t_1)
\end{array}
\right.}~ \eqno(2.10)
$$
Note that we can interpret a solution $y(t)$ of Eq. (2.7) on $[t_0,t_1)$ as a solution of a linear equation
$$
y' + H(t) y - R_1(t,\phi(t),\psi(t)) = 0
$$
on $[t_0,t_1)$, where $H(t) \equiv Q(t,\phi(t),\psi(t))y(t) + B(t,\phi(t),\psi(t)) + \frac{F(t,\phi(t),\psi(t))}{\phi(t)}, \linebreak R_1(t) \equiv~ R(t,\phi(t),\psi(t)) + \frac{G(t,\phi(t),\psi(t))}{\phi(t)}, \ph t\in [t_0,t_1)$. Then according to Cauchi formula we have
$$
y(t) = \exp\biggl\{-\il{t_0}{t}H(\tau) d \tau\biggr\}\biggl[y(t_0) + \il{t_0}{t}\exp\biggl\{\il{t_0}{\tau}H(s) d(s)\biggr\}R_1(\tau,\phi(\tau),\psi(\tau))d\tau\biggr],  \eqno (2.11)
$$
$t\in [t_0,t_1)$. By analogy for a solution $z(t)$ of Eq. (2.9) we get
$$
z(t) = \exp\biggl\{-\il{t_0}{t}L(\tau) d \tau\biggr\}\biggl[z(t_0) + \il{t_0}{t}\exp\biggl\{\il{t_0}{\tau}L(s) d(s)\biggr\}Q_1(\tau,\phi(\tau),\psi(\tau))d\tau\biggr],  \eqno (2.12)
$$
$t\in [t_0,t_1),$ where $L(t) \equiv R(t,\phi(t),\psi(t)) z(t) - B(t,\phi(t),\psi(t)) + \frac{G(t,\phi(t),\psi(t))}{\psi(t)}, \linebreak Q_1(t) \equiv~ Q(t,\phi(t),\psi(t)) + \frac{F(t,\phi(t),\psi(t))}{\psi(t)},  \ph t\in [t_0,t_1)$

{\bf Lemma 2.1.} {\it Let $(\phi(t),\psi(t))$ be a solution of the system (1.2) on $[t_0,t_1)$. If \linebreak $Q(t,\phi(t),\psi(t)) \ge~ 0, \ph R(t,\phi(t),\psi(t)) \ge 0, \ph
  t\in [t_0,t_1)$, then for every $\gamma \ge 0$ Eq. (2.7) (Eq. (2.9)) has a non negative solution $y(t) \ph (z(t))$ on $[t_0,t_1)$.}

Proof. Let us prove the existence of $y(t)$. Consider the Riccati equation
$$
y' + Q(t,\phi(t),\psi(t)) y^2 + B(t,\phi(t),\psi(t)) y = 0, \phh t\in [t_0,t_1) \eqno (2.13)
$$
Obviously $y_2(t) \equiv 0$ is a solution of this equation on $[t_0,t_1)$. Then since $Q(t,\phi(t),\psi(t)) \ge~ 0, \linebreak R(t,\phi(t),\psi(t)) \ge 0, \ph t\in [t_0,t_1)$, using Theorem 2.1 to the pair of equations (2.7) and (2.13) we conclude that for every $\gamma \ge 0$ Eq. (2.7) has a non negative solution $y(t)$ on $[t_0,t_1)$ with $y(t_0) = \gamma$.
The proof of existence of $z(t)$ can be made by analogy using Eq.~ (2.9) instead of Eq. (2.7). The lemma is proved.

{\bf Definition 2.1.} {\it An interval $[t_0,T)$ is called the maximum existence interval for a solution
$(\phi(t),\psi(t))$ of the system (1.1), if $(\phi(t),\psi(t))$ satisfies (1.1) almost everywhere  on $[t_0,T)$ and cannot be continued to the right from $T$ as a solution of Eq. (1.1).}

{\bf Lemma 2.2.} {\it Assume $P(t,u,v), \ph Q(t,u,v), \ph R(t,u,v), \ph S(t,u,v), \ph F(t,u,v)$ and \linebreak $G(t,u,v)$ are locally integrable and locally bounded  by $t$  and continuous by $u$ and $v$ on $[t_0,+\infty)\times \mathbb{R}\times \mathbb{R}.$ If $(\phi(t), \psi(t))$ is a bounded solution of the system (1.1) on $[t_0,T)$ and $T < +\infty$, then $[t_0,T)$ cannot be the maximum existence interval for $(\phi(t), \psi(t))$.}

Proof. To prove the lemma it is enough to show that  $(\phi(t), \psi(t))$ has a finite limit when $t \to T -0$. By (1.1) we have
$$
\phi(t) = \phi(t_0) + \il{t_0}{t}A_{\phi,\psi}(\tau) d \tau, \phh t \in [t_0,T), \eqno (2.14)
$$
where $A_{\phi,\psi}(\tau) \equiv P(\tau,\phi(\tau),\psi(\tau))\phi(\tau) + Q(\tau,\phi(\tau),\psi(\tau))\psi(\tau) +  F(\tau,\phi(\tau),\psi(\tau)), \ph \tau \in [t_0,T).$ Since $(\phi(t),\psi(t))$ is bounded on $[t_0,T)$ and $P(t,u,v), \ph Q(t,u,v), \ph R(t,u,v)$ are continuous by $u, \ph v$ on  $\mathbb{R}$ functions we have
$$
|A_{\phi,\psi}(t)| \le M, \phh t \in [t_0,T)    \eqno (2.15)
$$
for some $M \ge 0$. Let $\{t_k\}_{k=1}^{+\infty}$ be a sequence such that $t_k < T, \ph k =1, 2, ...$ and $\lim\limits_{k \to +\infty} t_k = T$. Then by (2.15) for any integers $m, n \ge 1$ we have $\biggl|\il{t_m}{t_n}A_{\phi,\psi}(\tau) d \tau\biggr| \le M |t_m - t_n|$. Hence by (2.14) $|\phi(t_m) - \phi(t_n)| \le M|t_m - t_n|,$ i. e. the sequence $\{\phi(t_k)\}_{k=1}^{+\infty}$ is fundamental, which proves that $\phi(t)$ has a finite limit when $t \to T - 0$. By analogy it can be proved the existence of the finite limit $\lim\limits_{t \to T - 0}\psi(t)$. The lemma is proved.

{\bf 3. Global solvability, stability and oscillation criteria}.

{\bf Definition 3.1.} {\it A solution $(\phi_0(t),\psi_0(t))$ of the system (1.1) ((1.2)) is called Lyapunov (asymptotically) stable if for every $\varepsilon > 0$ there exists $\delta > 0$ such that for every solution $(\phi(t), \psi(t))$ of the system (1.1) ((1.2)) with $|\phi(t_0) - \phi_0(t_0)| \le \delta, \ph |\psi(t_0) - \psi_0(t_0)| \le \delta$ the inequalities
$$
|\phi(t) - \phi_0(t)| \le \varepsilon, \phh |\psi(t) - \psi_0(t)| \le \varepsilon, \phh t \ge t_0
$$
are valid
(the relations
$$
\liml{t \to +\infty}(\phi(t) - \phi_0(t)) = \liml{t \to +\infty}(\psi(t) - \psi_0(t)) = 0
$$
are valid).
}

{\bf Definition 3.2.} {\it A solution $(\phi_0(t),\psi_0(t))$ of the system (1.1) ((1.2)) is called global asymptotically stable if for every its solution $(\phi(t), \psi(t))$ the functions $\phi(t) - \phi_0(t)$ and $\psi(t) - \psi_0(t)$ vanish at $+\infty$.}

Let $P_0(t), \ph Q_0(t), \ph R_0(t), \ph S_0(t), \ph F_0(t)$ and $G_0(t)$ be real-valued locally integrable functions on $[t_0,+\infty)$. Consider the linear system of ordinary differential equations
$$
\sist{\phi' = P_0(t) \phi + Q_0(t) \psi + F_0(t),}{\psi' = R_0(t)\phi + S_0(t) \psi + G_0(t), \ph t \ge t_0} \eqno (3.1)
$$
and the corresponding homogeneous one
$$
\sist{\phi' = P_0(t) \phi + Q_0(t) \psi,}{\psi' = R_0(t)\phi + S_0(t) \psi, \ph t \ge t_0} \eqno (3.2)
$$

{\bf Theorem 3.1.} {\it Assume $P(t,u,v), \ph Q(t,u,v), \ph R(t,u,v), \ph S(t,u,v), \ph  F(t,u,v)$ and \linebreak $G(t,u,v)$ are locally integrable and locally bounded  by $t$  and continuous by $u$ and $v$ on $[t_0,+\infty)\times \mathbb{R}\times \mathbb{R}.$ Moreover assume $P(t,u,v)\le P_0(t), \ph |Q(t,u,v)| \le Q_0(t), \linebreak |R(t,u,v)| \le~ R_0(t), \ph  S(t,u,v) \le S_0(t), \ph |F(t,u,v)| \le F_0(t), \ph  |G(t,u,v)| \le G_0(t), \linebreak t\ge~ t_0, u,v \in \mathbb{R}.$ Then the following statements are valid.

\noindent
I) For every $\alpha, \beta \in \mathbb{R}$  the solution $(\phi(t), \psi(t))$  of the system (1.1) with $\phi(t_0)=\alpha, \linebreak \psi(t_0) =~\beta$ exists on $]t_0,+\infty)$.

\noindent
II) If a  solution $(\phi_0(t), \psi_0(t))$ of the system (3.1) with $\phi_0(t_0) > 0, \ph \psi_0(t_0) > 0$ is bounden on $[t_0,+\infty)$ (vanish at $+\infty$), then every solution $(\phi(t), \psi(t))$ of the system (1.1) with $|\phi(t_0)| \le \phi_0(t_0), \ph |\psi(t_0)| \le \psi_0(t_0)$ is also bounded on $[t_0,+\infty)$ (vanish at $+\infty$).

\noindent
III) If the system (3.2) is Lyapunov (asymptotically) stable then the null solution $(0, 0)$ of the system (1.2) is Lyapunov (global asymptotically) stable.

}

Proof. Let us prove the statement I). Let $(\phi_0(t), \psi_0(t))$ be a solution of the system (1.1) with $\phi_0(t_0) = \alpha, \ph \psi_0(t_0) =~ \beta.$ Suppose this solutions does not exist on $[t_0,+\infty)$. Let then $[t_0,T)$ be the maximum existence interval for this solution. In virtue of (2.3) and (2.4) $\phi_0(t)$ and $\psi_0(t)$ are solutions of the Volterra equations
$$
\phi(t) = \stackrel{o}{v}_1(t) + \il{t_0}{t} \stackrel{o}{K}_1(t,\zeta)\phi(\zeta) d \zeta, \phh t\in [t_0,T),
$$
$$
\psi(t) = \stackrel{o}{v}_2(t) + \il{t_0}{t} \stackrel{o}{K}_2(t,\zeta)\psi(\zeta) d \zeta, \phh t\in [t_0,T)
$$
respectively, where
$$
\stackrel{o}{v}_1(t)\equiv \phi_0(t_0)\exp\biggl\{\il{t_0}{t}P(\tau,\phi_0(\tau),\psi_0(\tau))d\tau\biggr\} + \ph \psi_0(t_0)\il{t_0}{t}\exp\biggl\{\il{\tau}{t}P(s,\phi_0(s),\psi_0(s))ds +
 $$
 $$
 +\il{t_0}{\tau}S(s,\phi_0(s),\psi_0(s))d s\biggr\} d\tau + \il{t_0}{t}\exp\biggl\{\il{\tau}{t}P(s,\phi_0(s),\psi_0(s))d s\biggr\}\times \phantom{aaaaaaaaaaaaaaaaaaaaaaaaa}
$$
$$
\phantom{aaa}\times\biggl[F(\tau, \phi_0(\tau),\psi_0(\tau)) + Q(\tau,\phi_0(\tau),\psi_0(\tau))\il{t_0}{t}G(s,\phi_0(s),\psi_0(s))d s\biggr]d\tau
$$
$$
\stackrel{o}{v}_2(t)\equiv \psi_0(t_0)\exp\biggl\{\il{t_0}{t}S(\tau,\phi_0(\tau),\psi_0(\tau))d\tau\biggr\} + \ph \phi_0(t_0)\il{t_0}{t}\exp\biggl\{\il{\tau}{t}S(s,\phi_0(s),\psi_0(s))ds +
 $$
 $$
 +\il{t_0}{\tau}P(s,\phi_0(s),\psi_0(s))d s\biggr\} d\tau + \il{t_0}{t}\exp\biggl\{\il{\tau}{t}S(s,\phi_0(s),\psi_0(s))d s\biggr\}\times \phantom{aaaaaaaaaaaaaaaaaaaaaaaaa}
$$
$$
\phantom{aaa}\times\biggl[G(\tau, \phi_0(\tau),\psi_0(\tau)) + R(\tau,\phi_0(\tau),\psi_0(\tau))\il{t_0}{t}F(s,\phi_0(s),\psi_0(s))d s\biggr]d\tau
$$
$$
\stackrel{o}{K}_1(t,\zeta) \equiv R(t,\phi_0(t),\psi_0(t))\il{\zeta}{t}\exp\biggl\{\il{\tau}{t}P(s,\phi_0(s),\psi_0(s))d s + \il{\zeta}{\tau} S(s,\phi_0(s),\psi_0(s))d s\biggr\}\times
$$
$$
\phantom{aaaaaaaaaaaaaaaaaaaaaaaaaaaaaaaaaaaaaaaaaaaaaaaaaaaa}\times Q(\tau,\phi_0(\tau),\psi_0(\tau)) d \tau,
$$
$$
\stackrel{o}{K}_2(t,\zeta) \equiv Q(t,\phi_0(t),\psi_0(t))\il{\zeta}{t}\exp\biggl\{\il{\tau}{t}S(s,\phi_0(s),\psi_0(s))d s + \il{\zeta}{\tau} P(s,\phi_0(s),\psi_0(s))d s\biggr\}\times
$$
$$
\phantom{aaaaaaaaaaaaaaaaaaaaaaaaaaaaaaaaaaaaaaaaaaaaaaaaaaaa}\times R(\tau,\phi_0(\tau),\psi_0(\tau)) d \tau,
$$
It follows from the conditions of the theorem that
$$
|\stackrel{o}{v}_j(t)| \le \stackrel{oo}{v}_j(t), \ph |\stackrel{o}{K}_j(t,\zeta)| \le \stackrel{oo}{K}_j(t,\zeta), \ph j=1,2, \ph t,\zeta \in [t_0,T), \ph \zeta \le t, \eqno (3.3)
$$
where
$$
\stackrel{oo}{v}_1(t)\equiv |\phi_0(t_0)|\exp\biggl\{\il{t_0}{t}P_0(\tau)d\tau\biggr\} +  |\psi_0(t_0)|\il{t_0}{t}\exp\biggl\{\il{\tau}{t}P_0(s)ds
 +\il{t_0}{\tau}S_0(s)d s\biggr\} d\tau + \phantom{aaaaaaaaaa}
$$
$$
 \phantom{aaaaaaaaaaaaaaaaaaaaaaa} +\il{t_0}{t}\exp\biggl\{\il{\tau}{t}P_0(s)d s\biggr\}\biggl[F_0(\tau) + Q_0(\tau)\il{t_0}{t}G_0(s)d s\biggr]d\tau
$$
$$
\stackrel{oo}{v}_2(t)\equiv |\psi_0(t_0)|\exp\biggl\{\il{t_0}{t}S_0(\tau)d\tau\biggr\} +  |\phi_0(t_0)|\il{t_0}{t}\exp\biggl\{\il{\tau}{t}S_0(s)ds
 +\il{t_0}{\tau}P_0(s)d s\biggr\} d\tau + \phantom{aaaaaaaaaa}
$$
$$
 \phantom{aaaaaaaaaaaaaaaaaaaaaaa}+ \il{t_0}{t}\exp\biggl\{\il{\tau}{t}S_0(s)d s\biggr\}\biggl[G_0(\tau) + R_0(\tau)\il{t_0}{t}F_0(s)d s\biggr]d\tau
$$
$$
\stackrel{oo}{K}_1(t,\zeta) \equiv R_0(\zeta)\il{\zeta}{t}\exp\biggl\{\il{\tau}{t}P_0(s)d s + \il{\zeta}{\tau} S_0(s)d s\biggr\} Q_0(\tau) d \tau, \phantom{aaaaaaaaaaaaaaaaaaaaaaaa}
$$
$$
\stackrel{oo}{K}_2(t,\zeta) \equiv Q_0(\zeta)\il{\zeta}{t}\exp\biggl\{\il{\tau}{t}S_0(s)d s + \il{\zeta}{\tau} P_0(s)d s\biggr\} R_0(\tau) d \tau, \ph t,\zeta \in [t_0,T), \ph \zeta \le t.
$$
Let $\phi_{oo}(t)$ and  $\psi_{oo}(t)$ be solutions of the Volterra equations
$$
\phi(t) = \stackrel{oo}{v}_1(t) + \il{t_0}{t}\stackrel{oo}{K}_1(t,\zeta) \phi(\zeta) d \zeta, \phh t\in [t_0,T), \eqno (3.4)
$$
$$
\psi(t) = \stackrel{oo}{v}_2(t) + \il{t_0}{t}\stackrel{oo}{K}_2(t,\zeta) \psi(\zeta) d \zeta, \phh t\in [t_0,T), \eqno (3.5)
$$
respectively. Show that  $\phi_{oo}(t)$ and  $\psi_{oo}(t)$ are bounded functions on $[t_0,T)$. By (3.4) and (3.5)  for $\phi_{oo}(t)$ and  $\psi_{oo}(t)$ we have the representations via series expansion
$$
\phi_{oo}(t) = \stackrel{oo}{v}_1(t) + \il{t_0}{t}\stackrel{oo}{K}_1(t,\zeta) \stackrel{oo}{v}_1(\zeta) d\zeta  + \il{t_0}{\zeta}\stackrel{oo}{K}_1(t,\zeta) d \zeta \il{t_0}{\zeta}\stackrel{oo}{K}_1(t,\xi)\stackrel{oo}{v}_1(\xi) d\xi + ... , \eqno (3.6)
$$
$$
\psi_{oo}(t) = \stackrel{oo}{v}_2(t) + \il{t_0}{t}\stackrel{oo}{K}_2(t,\zeta) \stackrel{oo}{v}_2(\zeta) d\zeta  + \il{t_0}{\zeta}\stackrel{oo}{K}_2(t,\zeta) d \zeta \il{t_0}{\zeta}\stackrel{oo}{K}_2(t,\xi)\stackrel{oo}{v}_2(\xi) d\xi + ... , \eqno (3.7)
$$
respectively. From the definitions of $\stackrel{oo}{v}_j(t), \ph \stackrel{oo}{K}_j(t,\zeta), \ph j=1,2$ is seen that they are bounded functions for $t,\zeta \in [t_0,T), \ph \zeta \le t$. Then from (3.6) and  (3.7) we obtain respectively:
$$
|\phi_{oo}(t)| \le m_1\biggl(1 + M_1(T - t_0) + \frac{[M_1(T - t)]^2}{2!} + ... \biggr) = m_1\exp\Bigl\{M_1(T-t_0)\Bigr\}, \eqno (3.8)
$$
$$
|\psi_{oo}(t)| \le m_2\biggl(1 + M_1(T - t_0) + \frac{[M_2(T - t)]^2}{2!} + ... \biggr) = m_2\exp\Bigl\{M_2(T-t_0)\Bigr\}, \eqno (3.9)
$$
where $m_j\equiv \sup\limits_{t\in[t_0,T)}|\stackrel{oo}{v}_j(t)|, \ph M_j \equiv \sup\limits_{t,\zeta \in [t_0,T), \ph \zeta \le t}|\stackrel{oo}{K}_j(t,\zeta)|, \ph j=1,2.$ Hence $\phi_{oo}(t)$ and  $\psi_{oo}(t)$ are bounded functions on $[t_0,T)$. This together  with (3.3), (3.8) and (3.9) implies that $(\phi_0(t),\psi_0(t))$ is bounded on $[t_0,T)$. Then by virtue of Lemma 2.2  $[t_0,T)$ is not the maximum existence interval for  $(\phi_0(t),\psi_0(t))$, which contradicts our assumption. The obtained contradiction completes the proof of the statement I). The statement III) follows immediately from the statement II). Prove the statement II). Let $(\phi_0(t), \psi_0(t))$ be a solution of the system (3.1) with $\phi_0(t_0) > 0, \ph \psi_0(t_0) > 0$. Then $\phi_0(t)$ and $\psi_0(t)$ are solutions of the Volterra equations (3.4) and (3.5) respectively on $[t_0,+\infty)$. Let $(\phi(t), \psi(t))$ be a solution of the system (1.1) with
$$
|\phi(t_0)| \le \phi_0(t_0), \ph |\psi(t_0)| \le \psi_0(t_0) \eqno (3.10)
$$
(according to the statement I) $(\phi(t), \psi(t))$ exists on $[t_0,+\infty)$). In virtue of (2.3) and (2.4) for $\phi(t)$ and $\psi(t)$ the relations
$$
\phi(t)= v_1(t) + \il{t_0}{t}K_1(t,\zeta)\phi(\zeta) d \zeta, \phh t \ge t_0,
$$
$$
\psi(t)= v_2(t) + \il{t_0}{t}K_2(t,\zeta)\psi(\zeta) d \zeta, \phh t \ge t_0
$$
are valid respectively, where
$$
v_1(t)\equiv \phi(t_0)\exp\biggl\{\il{t_0}{t}P(\tau,\phi(\tau),\psi(\tau))d\tau\biggr\} + \ph \psi(t_0)\il{t_0}{t}\exp\biggl\{\il{\tau}{t}P(s,\phi(s),\psi(s))ds +
 $$
 $$
 +\il{t_0}{\tau}S(s,\phi(s),\psi(s))d s\biggr\} d\tau + \il{t_0}{t}\exp\biggl\{\il{\tau}{t}P(s,\phi(s),\psi(s))d s\biggr\}\times \phantom{aaaaaaaaaaaaaaaaaaaaaaaaa}
$$
$$
\phantom{aaa}\times\biggl[F(\tau, \phi(\tau),\psi(\tau)) + Q(\tau,\phi(\tau),\psi(\tau))\il{t_0}{t}G(s,\phi(s),\psi(s))d s\biggr]d\tau
$$
$$
v_2(t)\equiv \psi(t_0)\exp\biggl\{\il{t_0}{t}S(\tau,\phi(\tau),\psi(\tau))d\tau\biggr\} + \ph \phi(t_0)\il{t_0}{t}\exp\biggl\{\il{\tau}{t}S(s,\phi(s),\psi(s))ds +
 $$
 $$
 +\il{t_0}{\tau}P(s,\phi(s),\psi(s))d s\biggr\} d\tau + \il{t_0}{t}\exp\biggl\{\il{\tau}{t}S(s,\phi(s),\psi(s))d s\biggr\}\times \phantom{aaaaaaaaaaaaaaaaaaaaaaaaa}
$$
$$
\phantom{aaa}\times\biggl[G(\tau, \phi(\tau),\psi(\tau)) + R(\tau,\phi(\tau),\psi(\tau))\il{t_0}{t}F(s,\phi(s),\psi(s))d s\biggr]d\tau
$$
For $\phi_0(t)$ and $\psi_0(t)$ we have the relations
$$
\phi_0(t) = \stackrel{0}{u}_1(t) +\il{t_0}{t}\stackrel{0}{L}_1(t,\zeta)\phi_0(\zeta) d\zeta, \phh t \ge t_0,
$$
$$
\psi_0(t) = \stackrel{0}{u}_2(t) +\il{t_0}{t}\stackrel{0}{L}_2(t,\zeta)\psi_0(\zeta) d\zeta, \phh t \ge t_0,
$$
respectively, where
$$
\stackrel{0}{u}_1(t)\equiv \phi_0(t_0)\exp\biggl\{\il{t_0}{t}P_0(\tau)d\tau\biggr\} +  \psi_0(t_0)\il{t_0}{t}\exp\biggl\{\il{\tau}{t}P_0(s)ds
 +\il{t_0}{\tau}S_0(s)d s\biggr\} d\tau + \phantom{aaaaaaaaaa}
$$
$$
 \phantom{aaaaaaaaaaaaaaaaaaaaaaa} +\il{t_0}{t}\exp\biggl\{\il{\tau}{t}P_0(s)d s\biggr\}\biggl[F_0(\tau) + Q_0(\tau)\il{t_0}{t}G_0(s)d s\biggr]d\tau
$$
$$
\stackrel{0}{u}_2(t)\equiv \psi_0(t_0)\exp\biggl\{\il{t_0}{t}S_0(\tau)d\tau\biggr\} +  \phi_0(t_0)\il{t_0}{t}\exp\biggl\{\il{\tau}{t}S_0(s)ds
 +\il{t_0}{\tau}P_0(s)d s\biggr\} d\tau + \phantom{aaaaaaaaaa}
$$
$$
 \phantom{aaaaaaaaaaaaaaaaaaaaaaa}+ \il{t_0}{t}\exp\biggl\{\il{\tau}{t}S_0(s)d s\biggr\}\biggl[G_0(\tau) + R_0(\tau)\il{t_0}{t}F_0(s)d s\biggr]d\tau
$$
$$
\stackrel{0}{L}_1(t,\zeta) \equiv R_0(\zeta)\il{\zeta}{t}\exp\biggl\{\il{\tau}{t}P_0(s)d s + \il{\zeta}{\tau} S_0(s)d s\biggr\} Q_0(\tau) d \tau, \phantom{aaaaaaaaaaaaaaaaaaaaaaaa}
$$
$$
\stackrel{0}{L}_2(t,\zeta) \equiv Q_0(\zeta)\il{\zeta}{t}\exp\biggl\{\il{\tau}{t}S_0(s)d s + \il{\zeta}{\tau} P_0(s)d s\biggr\} R_0(\tau) d \tau, \ph t \ge \zeta \ge t_0.
$$
It follows from the conditions of the theorem that
$$
|v_j(t)| \le \stackrel{0}{u}_j(t), \ph |k_j(t,\zeta)| \le |\stackrel{0}{L}_j(t)|, \ph j=1,2, \ph t \ge \zeta \ge t_0.
$$
Then from (3.10) it follows that
$$
|\phi(t)| \le \phi_0(t), \phh |\psi(t)| \le \psi_0(t), \phh t \ge t_0. \eqno (3.11)
$$
Assume $(\phi_0(t), \psi_0(t))$ is bounded on $[t_0,+\infty)$ (vanish at $+\infty$). Then the relations (3.11) imply that $(\phi(t), \psi(t))$ is bounded on $[t_0,+\infty)$ (vanish at $+\infty$). The statement II) is proved. The proof of the theorem is complete.

{\bf Remark 3.1.} {\it From the proof of Theorem 3.1 is seen that it remains valid if we generalize its conditions assuming that
$P(t,u,v), \ph Q(t,u,v), \ph R(t,u,v), \ph S(t,u,v), \ph F(t,u,v)$ and  $G(t,u,v)$ are complex-valued
 locally integrable and locally bounded by $t$ and continuous by $u$ and $v$  functions on
$[t_0,+\infty) \times (\mathbb{C})\times  (\mathbb{C})$ and replace  $\alpha, \beta \in \mathbb{R}$ by $\alpha, \beta \in \mathbb{C}$ and  the conditions $P(t,u,v) \le P_0(t), \ph S(t,u,v) \le S_0(t)$ by  $Re \hskip 2pt P(t,u,v) \le P_0(t), \linebreak Re \hskip 2pt  S(t,u,v) \le~ S_0(t)$ .}

{\bf Example 3.1}. {\it  Eq. (1.3) is equivalent to the system
$$
\sist{\phi' = \psi,}{\psi' = -(1 + \beta cos t)\phi - \varepsilon (\phi^2 -1) \psi, \ph t\ge t_0.}
$$
Obviously for this system with $\varepsilon \ge 0$ the conditions of Theorem 3.1 are satisfied. Therefore for every $\alpha, \beta \in \mathbb{R}$ Eq. (1.3) has a solution $\phi(t)$ on $[t_0,+\infty)$ with $\phi(t_0) = \alpha, \ph \psi(t_0) = \beta$. By similar way can be discussed the applicability of Theorem 3.1 to the equations (1.4) - (1.14). The applicability of Theorem 3.1 to the
 system of equations (see [1], p. 257):
$$
\sist{\phi'= -(\phi^2 + \psi^2 -1) \phi + a \sin \hskip 2pt t \hskip 2pt\psi,}{\psi'= - a \sin \hskip 2pt t \hskip 2pt\phi - (\phi^2 + \psi^2 -1)\psi}
$$
is obvious.}

Let $B_1(t)$ and $B_2(t)$ be real-valued continuous functions on $[t_0,+\infty)$. For any $c_1 \ne 0$ and $c_2 \ne 0$ set:
$$
K(t,c_1,c_2) \equiv c_1\exp\biggl\{\il{t_0}{t}P_0(s)d s + \il{t_0}{t}Q_0(\tau)\biggr[\frac{c_2}{c_1}\exp\biggl\{-\il{t_0}{\tau}B_1(s) d s\biggr\} + \phantom{aaaaaaaaaaaaaaaaaa}
$$
$$
\phantom{aaaaaaaaaaaaaaaaaaaaaaaaaaaaaaaaaaaaaaaa}+\il{t_0}{\tau}\exp\biggl\{-\il{\zeta}{\tau}B_1(s) d s\biggr\} R_0(\zeta) d \zeta\biggr]d\tau\biggr\},
$$
$$
L(t,c_1,c_2) \equiv c_2\exp\biggl\{\il{t_0}{t}S_0(s)d s + \il{t_0}{t}R_0(\tau)\biggr[\frac{c_1}{c_2}\exp\biggl\{\il{t_0}{\tau}B_2(s) d s\biggr\} + \phantom{aaaaaaaaaaaaaaaaaa}
$$
$$
\phantom{aaaaaaaaaaaaaaaaaaaaaaaaaaaaaaaaa}+\il{t_0}{\tau}\exp\biggl\{\il{\zeta}{\tau}B_2(s) d s\biggr\} Q_0(\zeta) d \zeta\biggr]d\tau\biggr\}, \phh t\ge t_0.
$$

{\bf Theorem 3.2.}. {\it Let for some $c_1 > 0, \ph c_2> 0, \ph \varepsilon > 0$ and for every $t\ge t_0, \linebreak u~ \in~ (0, K(t,c_1,c_2) + \varepsilon], \ph v \in (0, L(t,c_1,c_2) + \varepsilon]$ the inequalities $P(t,u,v) \le P_0(t), \linebreak S(t,u,v) \le~ S_0(t), \ph 0\le Q(t,u,v) \le Q_0(t), \ph 0\le R(t,u,v) \le R_0(t), \ph B_1(t) \le B(t,u,v) \le B_2(t)$ are valid. Then  every  solution $(\phi(t),\psi(t))$ of the system (1.2)  with $\phi(t_0) = c_1, \linebreak \psi(t_0) =~ c_2$ exists  on $[t_0,+\infty)$  and
$$
0 < \phi(t) \le K(t,c_1,c_2), \phh 0< \psi(t) \le L(t,c_1,c_2), \phh t\ge t_0. \eqno (3.12)
$$
}

Proof. Let $(\phi_0(t),\psi_0(t))$ be a solution of the system (1.2) with $\phi_0(t_0)= c_1$ and $\psi_0(t_0)=~c_2.$ Suppose this solution is not continuable on $[t_0,+\infty)$. Let then $[t_0,T)$ be the maximum existence interval for $(\phi_0(t),\psi_0(t))$. Show that the inequalities (3.12) are satisfied for all $t\in [t_0,T)$. By (2.1) and (2.2) we have
$$
\phi_0(t) > 0, \ph \psi_0(t) > 0, \ph t\in [t_0,T). \eqno (3.13)
$$
Note that (3.12) is valid at least for $t = t_0$. Then we can set:
$$
\overline{t}_1\equiv \sup\{t\in [t_0,T): \phi_0(\tau) \le K(\tau,c_1,c_2), \ph \tau \in [t_0,t]\},
$$
$$
\overline{t}_2\equiv \sup\{t\in [t_0,T): \psi_0(\tau) \le L(\tau,c_1,c_2), \ph \tau \in [t_0,t]\},
$$
Assume $\overline{t}_1 \le \overline{t}_2$ (the proof in the case  $\overline{t}_2 \le \overline{t}_1$ by analogy) and assume (3.12) is not valid for all $t\in [t_0,T)$. Then $\overline{t}_1 < T$ and, therefore, there exists $t_2\in (t_0,T)$ such that
$$
\phi_0(t_2) > K(t_2,c_1,c_2), \eqno (3.14)
$$
$$
\phi_0(t) \le K(t,c_1,c_2) + \varepsilon, \ph t\in [t_0,t_2], \eqno (3.15)
$$
$$
\psi_0(t) \le L(t,c_1,c_2) + \varepsilon, \ph t\in [t_0,t_2]. \eqno (3.16)
$$
Consider the Riccati equations
$$
y' + Q(t,\phi_0(t),\psi_0(t)) y^2  + B(t,\phi_0(t),\psi_0(t)) y - R(t,\phi_0(t), \psi_0(t)) = 0, \eqno (3.17)
$$
$$
z' + R(t,\phi_0(t),\psi_0(t)) z^2  - B(t,\phi_0(t),\psi_0(t)) z - Q(t,\phi_0(t), \psi_0(t)) = 0, \eqno (3.18)
$$
By Lemma 2.1 from the conditions of the theorem and from (3.15), (3.16) it follows that for every $\gamma \ge 0,$ Eq. (3.17) (Eq. (3.18)) has a non negative solution $y_0(t) \ph (z_0(t))$ on $[t_0,t_2)$  with $y_0(t_0) = \gamma \ph (z_0(t_0) = \gamma)$. By (2.11) we have
$$
y_0(t)= y_0(t_0)\exp\biggl\{-\il{t_0}{t}H_0(\tau)d\tau\biggr\} + \il{t_0}{t}\exp\biggl\{-\il{\tau}{t}H_0(s) d s\biggr\}R(\tau,\phi_0(\tau),\psi_0(\tau)) d\tau,
$$
$t\in [t_0,t_2)$, where $H_0(t) \equiv Q(t,\phi_0(t),\psi_0(t)) y_(t) + B(t,\phi_0(t),\psi_0(t)), \ph t\in [t_0,t_2).$
Multiply both sides of this equality by $Q(t,\phi_0(t),\psi_0(t))$ and integrate from $t_0$ to $t$. We obtain
$$
\il{t_0}{t}Q(\tau,\phi_0(\tau),\psi_0(\tau)) y_0(\tau) d\tau = \il{t_0}{t}Q(\tau,\phi_0(\tau),\psi_0(\tau))\times \phantom{aaaaaaaaaaaaaaaaaaaaaaaaaaaaaaaaaaaaaaaaa}
$$
$$
\phantom{aaa}\biggl[y_0(t_0)\exp\biggl\{-\il{t_0}{\tau}H_0(s) d s\biggr\} + \il{t_0}{\tau}\exp\biggl\{-\il{\zeta}{\tau}H_0(s) d s\biggr\}R(\zeta,\phi_0(\zeta),\psi_0(\zeta)) d\zeta\biggr]d \tau,
$$
$t\in [t_0,t_2).$ By (2.8) this equality with conditions of the theorem implies
$$
\phi_0(t_2) \le K(t_2,c_1,c_2),
$$
which contradicts (3.14). The obtained contradiction implies
$$
0 < \phi(t) \le K(t,c_1,c_2), \phh 0< \psi(t) \le L(t,c_1,c_2), \phh t\in [t_0,T). \eqno (3.19)
$$
From here and from (3.13) it follows that $(\phi_0(t), \psi_0(t))$ is bounded on $[t_0,T)$. Then by Lemma 2.2  $[t_0,T)$ is not the maximum existence interval for $(\phi_0(t),\psi_0(t))$, which contradicts our assumption. The obtained contradiction together with (3.13) and  (3.19) completes the proof of the theorem.

{\bf Definition 3.3.} {\it A solution $(\phi(t), \psi(t))$ of the system (1.1) is called conditionally stable, if there exists an one dimensional manifold $J \ni (\phi(t_0), \psi(t_0))$ such that for every $\varepsilon > 0$ there exists $\delta >0$ such that each solution $(\widetilde{\phi}(t), \widetilde{\psi}(t))$ of the system (1.1) with $(\widetilde{\phi}(t_0), \widetilde{\psi}(t_0)) \in J$ and $|\widetilde{\phi}(t_0) - \phi(t_0)| + |\widetilde{\psi}(t_0) - \psi(t_0)|  < \delta$ exists on $[t_0,+\infty)$ and satisfies to the inequality   $|\widetilde{\phi}(t) - \phi(t| + |\widetilde{\psi}(t) - \psi(t)|  < \varepsilon, \ph t \ge t_0.$.}

{\bf Corollary 3.1.} {\it Let for some $c_1 > 0, \ph c_2 > 0, \ph \varepsilon > 0$ and for every $t \ge t_0, \ph \pm u \in (0, K(t,c_1, c_2) + \varepsilon], \ph \pm v \in (0, L(t, c_1,c_2) + \varepsilon]$ the inequalities $P(t, \pm u, \pm v) \le P_0(t), \linebreak S(t, \pm u, \pm v) \le S_0(t), \ph  0 \le Q(t, \pm u, \pm v) \le Q_0(t), \ph 0 \le  R(t, \pm u, \pm v) \le R_0(t), \ph B_1(t) \le B(t, \pm u, \pm v) \le B_2(t)$ are valid (the expression $(t, \pm u, \pm v)$ means $(t, u, v)$ and $(t, -u, -v)$) and let $\lim\limits_{u \to 0, u\ne 0} P(t,u,v) u = 0, \ph  \lim\limits_{v \to 0, v\ne 0} S(t,u,v) v = 0$.  Then every solution $(\phi(t), \psi(t))$ of the system (1.2) with $\phi(t_0) = \lambda c_1, \ph \psi(t_0) = \lambda c_2, \ph |\lambda| \le~ 1$ exists on $[t_0, +\infty)$ and
$$
0 <  |\phi(t)| \le \lambda K(t, c_1, c_2), \ph 0 < |\psi(t)| \le \lambda L(t,c_1,c_2), \phh t \ge t_0. \eqno (3.20)
$$
Moreover if $K(t, c_1, c_2)$ and $L(t, c_1, c_2)$ are bounded on $[t_0,+\infty)$ ($c_1$ and $c_2$ are fixed) then the null solution $(0,0)$  of the system (1.2) is conditionally stable.
}

Proof.
The substitution $\phi \rightarrow - \phi, \ph \psi \rightarrow - \psi$ in the system (1.2) reduces it in the following one
$$
\sist{\phi'=P(t,-\phi,-\psi)\phi + Q(t,-\phi,-\psi)\psi,}{\psi' = R(t,-\phi,-\psi)\phi + S(t,-\phi,-\psi)\psi, \ph t \ge t_0.} \eqno (3.21)
$$
Moreover it is obvious that
$$
K(t,\lambda c_1, \lambda c_2) = \lambda K(t, c_1, c_2), \phh L(t,\lambda c_1, \lambda c_2) = \lambda L(t, c_1, c_2), \ph t \ge t_0, \ph \lambda \in \mathbb{R}
$$
By Theorem 3.2 it follows from here and from $(3.21)$ that under the conditions of the corollary every solution $(\phi(t), \psi(t))$ of the system (1.2) with $\phi(t_0) = \lambda c_1, \ph \psi(t_0) = \lambda c_2, \linebreak |\lambda| \le 1$ exists on $[t_0,+\infty)$ and the estimates $(3.20)$ are valid. If as a manifold, indicated in Definition 3.3, we take $J \equiv \{(\lambda c_1, \lambda c_2)\}$ then from $(3.20)$ and from the boundedness of $K(t, c_1, c_2)$ and $L(t, c_1, c_2)$ it follows the conditionally stable of the null solution $(0,0)$ of the system (1.2). The corollary is proved.

Let $\alpha(t), \ph \beta(t), \ph P_\delta(t), \ph Q_\delta(t), \ph R_\delta(t), \ph S_\delta(t), \ph B_{1,\delta}(t), \ph B_{2,\delta}(t)$ be real-valued locally \linebreak integrable functions on $[t_0,+\infty), \ph \delta > 0$. For any $c_1  \ne 0, \ph c_2 \ne 0$ set:

$$
K_\delta(t,c_1,c_2) \equiv c_1\exp\biggl\{\il{t_0}{t}P_\delta(s)d s + \il{t_0}{t}Q_\delta(\tau)\biggr[\frac{c_2}{c_1}\exp\biggl\{-\il{t_0}{\tau}B_{1,\delta}(s) d s\biggr\} + \phantom{aaaaaaaaaaaaaaaaaa}
$$
$$
\phantom{aaaaaaaaaaaaaaaaaaaaaaaaaaaaaaaaaaaaaaaa}+\il{t_0}{\tau}\exp\biggl\{-\il{\zeta}{\tau}B_{1,\delta}(s) d s\biggr\} R_\delta(\zeta) d \zeta\biggr]d\tau\biggr\},
$$
$$
L_\delta(t,c_1,c_2) \equiv c_2\exp\biggl\{\il{t_0}{t}S_\delta(s)d s + \il{t_0}{t}R_\delta(\tau)\biggr[\frac{c_1}{c_2}\exp\biggl\{\il{t_0}{\tau}B_{2,\delta}(s) d s\biggr\} + \phantom{aaaaaaaaaaaaaaaaaa}
$$
$$
\phantom{aaaaaaaaaaaaaaaaaaaaaaaaaaaaaaaaa}+\il{t_0}{\tau}\exp\biggl\{\il{\zeta}{\tau}B_{2,\delta}(s) d s\biggr\} Q_\delta(\zeta) d \zeta\biggr]d\tau\biggr\}, \phh t\ge t_0.
$$

{\bf Theorem 3.3.} {\it Let for some $c_1 > 0, \ph c_2 > 0, \ph \varepsilon > 0$ and for every $\delta \in (0, \varepsilon), \ph t \ge t_0 , \ph u \in [\delta, K_\delta(t, c_1, c_2)+ \varepsilon], \ph v \in [\delta, L_\delta(t, c_1, c_2) + \varepsilon ]$ the inequalities
$$
\alpha(t) \le P(t,u,v) \le P_\delta(t), \ph \beta(t) \le S(t,u,v) \le S_\delta(t)  \eqno (3.22)
$$
$0 \le Q(t,u,v) \le Q_\delta(t), \ph 0 \le R(t,u,v) \le R_\delta(t), \ph B_{1,\delta}(t) \le B(t, u,v) \le B_{2,\delta}(t)$ are valid.  Then every solution $(\phi(t), \psi(t))$ of the system (1.2) with $\phi(t_0) = c_1, \ph \psi(t_0) = c_2$ exists on $[t_0, +\infty)$ and
$$
\phi(t) > 0, \phh \psi(t) > 0, \phh t \ge t_0. \eqno (3.23)
$$
}

Proof. Let $(\phi_0(t), \psi_0(t))$ be a solution of the system (1.2) with $\phi_0(t_0) = c_1, \ph \psi_0(t_0) = c_2$ and let $[t_0, T)$ be its maximum existence interval. We must show that
$$
T = + \infty. \eqno (3.24)
$$
By (2.1) and (2.2) from the conditions $(3.22)$ it follows that
$$
\phi_0(t) \ge c_1\exp\biggl\{\il{t_0}{t} \alpha(\tau) d \tau\biggr\}, \phh
\psi_0(t) \ge c_2\exp\biggl\{\il{t_0}{t} \beta(\tau) d \tau\biggr\}, \phh t \in [t_0,T);
$$
Suppose $T<+\infty.$ Then since $\alpha(t)$ and $\beta(t)$ are locally integrable  from the last inequalities it follows that
$$
\phi_0(t) \ge \delta_0, \phh \psi_0(t) \ge \delta_0, \phh t \in [t_0, T)
$$
for some $\delta_0 > 0$. Show that
$$
0 < \phi_0(t) \le K_{\delta_0}(t, c_1, c_2), \phh 0 < \psi_0(t) \le L_{\delta_0}(t, c_1, c_2), \phh t \in [t_0,T). \eqno(3.25)
$$
From the conditions of the theorem is seen that these inequalities are satisfied at least for $t = t_0$. Then we can set:
$$
\overline{t}_1\equiv \sup\{t\in [t_0,T): \phi_0(\tau) \le K_{\delta_0}(\tau,c_1,c_2), \ph \tau \in [t_0,t]\},
$$
$$
\overline{t}_2\equiv \sup\{t\in [t_0,T): \psi_0(\tau) \le L_{\delta_0}(\tau,c_1,c_2), \ph \tau \in [t_0,t]\},
$$
Assume $\overline{t}_1 < \overline{t}_2$ (the proof in the case $\overline{t}_2 \le \overline{t}_1$ by analogy) and assume $(\delta)$ is not valid Then $\overline{t}_t < T$ and there exists $t_2 \in (t_0, T)$ such that
$$
\phi_0(t_2) > K_{\delta_0}(t_2, c_1, c_2), \eqno (3.26)
$$
$$
\delta_0 \le \phi_0(t) \le K_{\delta_0}(t,c_1,c_2) + \varepsilon, \ph t \in [t_0, t_2],
$$
$$
\delta_0 \le \psi_0(t) \le L_{\delta_0}(t,c_1,c_2) + \varepsilon, \ph t \in [t_0, t_2].
$$
Further similar to the proof of Theorem 3.2 from here we come to the conclusion that
$$
\phi_0(t_2) \le K_{\delta_0}(t,c_1, c_2),
$$
which contradicts (3.26). The obtained contradiction  proves (3.25). It follows from (3.25)  that $(\phi_0(t), \psi_0(t))$ is bounded on $[t_0,T).$. Then by Lemma 2.2 (3.24) is valid. from (3.24)   and from (3.25)  it follows (3.23). The theorem is proved.

In contrast to the conditions of Theorem 3.2, the conditions of Theorem 3.3 admit that the functions $ Q $ and $ R $ have singularities for $ u = 0 $ and $ v = 0 $, which, at part, is illustrated in the following example.

{\bf Example 3.2.} {\it Consider the Emden-Fowler's equation (see [8], p. 171)
$$
(t^\rho\phi')' - t^\sigma \phi^n = 0, \phh t \ge t_0 > 0,  \phh \rho,\sigma, n \in \mathbb{R}. \eqno  (3.27)
$$
This equation is equivalent to the following system
$$
\sist{\phi'= t^{-\rho} \psi,}{\psi' = (t^\sigma \phi^{n-1})\phi, \ph t\ge t_0.} \eqno (3.28)
$$
Here $P(t,u,v) = S(t,u,v) \equiv 0, \ph Q(t,u,v) = t^{-\rho}, \ph R(t,u,v) = t^\sigma u^{n-1}$. Consider the case $\rho \ne 1, \ph \sigma \ne -1, \ph 2 + \sigma - \rho \ne 0.$ At first we consider the superlinear case $ n > 1$. Let us take $P_0(t) = S_0(t)\equiv 0, \ph Q_0(t)= t^{-\rho}, \ph R_0(t) = t^\sigma, \ph t\ge t_0 > 0.$ Then it is not difficult to show that
$$
L(t,c_1,c_2) = c_2\exp\biggl\{\biggl(\frac{c_1}{c_2} - \frac{t_0^{1-\rho}}{1-\rho}\biggr)\frac{t^{\sigma + 1} - t_0^{\sigma + 1}}{\sigma + 1} + \frac{t^{2 + \sigma - \rho} - t_0^{2 + \sigma - \rho}}{(1 -\rho)(2 + \sigma - \rho)}\biggr\},  \eqno  (3.29)
$$
$t\ge t_0 > 0.$
Obviously the inequality $0\le Q(t,u,v) \le Q_0(t), \ph t\ge t_0 > 0$ for the system  (3.28) holds for all $u,v \in\mathbb{R}$. Then on the basis of  (3.29) it is easy to verify that the conditions of Theorem 3.2 for the system  (3.28) are satisfied provided:

\noindent
1. $\rho < 1, \ph \sigma < -1, \ph 2 + \sigma - \rho  < 0, \ph c_1 > 0, \ph 0 < c_2 < 1, \ph \frac{c_1}{c_2} \le \frac{t_0^{1-\rho}}{1 - \rho},$

\noindent
or

\noindent
2. $\rho <1, \ph \sigma < -1, \ph 2 + \sigma - \rho < 0, \ph c_1 > 0, \ph 0 < c_1 < \exp\biggl\{\frac{t_0^{\sigma + 1}}{\sigma + 1}\biggl\}, \ph \frac{c_1}{c_2} > \frac{t_0^{1-\rho}}{1 - \rho}. $

\noindent
In the case $\rho = 1, \ph \sigma \ne -1$ we have
$$
L(t,c_1,c_2) = c_2\exp\biggl\{\frac{c_1}{c_2}\frac{t^{\sigma + 1} - t_0^{\sigma + 1}}{\sigma + 1} + \il{t_0}{t}\tau^\sigma \ln \frac{\tau}{t_0}d \tau\biggr\}.
$$
Hence the conditions of Theorem 3.2 for the system  (3.28) are satisfied provided

\noindent
3. $\rho =1, \ph \sigma < -1, \ph c_1 > 0, \ph 0 < c_2 < \exp\biggl\{\frac{c_1}{c_2}\frac{t_0^{\sigma + 1}}{\sigma + 1} - \ilp{t_0}\tau^\sigma \ln \frac{\tau}{t_0}d \tau\biggr\}.$
In the sublinear ($0 < n <1$) and extra sublinear ($n < 0$) cases (the cases $n=0$ and $n=1$ are trivial) we can take $\alpha(t) = \beta(t) = P_\delta(t) = S_\delta(t) \equiv 0, \ph Q_\delta(t) = t^{-\rho}, \ph R_\delta(t) = t^\sigma \delta^{n-1}$ and apply Theorem 3.3.
}

{\bf Remark 3.2.} {\it
The case $\rho >1, \ph n \ne 1$ is reducible to the considered case $\rho < 1$ (namely to its simpler case $\rho = 0$) by the transformation (see [8], p. 143)
$$
s = \frac{t^{\rho -1}}{\rho -1}, \phh \phi = (\rho -1)^{(\rho - \sigma -2)/ [(\rho - 1)(n-1)]}\hskip 2pt \frac{\psi}{s} \eqno (3.30)
$$
which reduces Eq. (3.22) to the following one
$$
\psi'' - s^{\sigma_1} \psi^n  = 0, \eqno (3.31)
$$
where $\sigma_1 \equiv \frac{\sigma + \rho}{\rho -1} - (n+3)$. Note that in the case $\rho <1, \ph n \ne 1$ Eq. (3.27) is also reducible to its simpler form by the another transformation (see [8], p. 144)
$$
\psi'' - s^{\sigma_2} \psi^n  = 0, \eqno (3.32)
$$
where $\sigma_2 \equiv \frac{\sigma + \rho}{1 - \rho}$. It was established in [8, p. 144] that for each pair of certain parameters $\rho$ and $n$ ($(\sigma + \rho)(\sigma +n -1) > 0$) the equation
$$
\phi'' - t^\sigma \phi^n = 0
$$
has a global solution of the form
$$
\phi(t) \equiv c(\sigma,n) t^{w(\sigma,n)}
$$
where $c(\sigma,n) \equiv \biggl[\frac{(\sigma = 2)(\sigma + n +1)}{(n -1)^2}\biggr]^{1/(n -1)}, \ph w(\sigma,n) \equiv -\frac{\sigma +2}{n-1}$.
Due to $(3.30)-(3.32)$ this means, that it is established the existence of a single global solution of Eq. (3.27) foe each ordered  triple $(\rho, \sigma, n)$ of certain parameters   $\rho, \ph \sigma, \ph n$. Whereas, in Example 3.2 has been \linebreak established the existence of two-parameter families of global solutions  of Eq. (3.27) for certain  parameters   $\rho, \ph \sigma, \ph n$.
}

{\bf Remark 3.3.} {\it The Wintner's theorem (see [9], pp. 29,30, Theorem 5.1) is not applicable neither to the van der Pol's nor to the Emden-Fowler's equations. Therefore on the basis of examples 3.1 and 3.2 we conclude that Theorem 3.1 and Theorem 3.2 are not consequences  of the Wintner's theorem.}

{\bf Theorem 3.4.} {\it Let the following conditions be satisfied.

\noindent
1) $0 \le Q(t,u,v) \le Q_0(t), \ph t \ge t_0, \ph u > 0;$

\noindent
2) $P(t,u,v) \le P_0(t), \ph t \ge t_0, \ph u > 0;$

\noindent
3) $0 \le F(t,u,v) \le F_0(t), \ph t \ge t_0, \ph u > 0;$

\noindent
4) $B(t,u,v) + \frac{F(t,u,v)}{u} \le L_0(t), \ph t \ge t_0, \ph u > 0;$

\noindent
5) $0 \le R(t,u,v) + \frac{G(t,u,v)}{u} \le M_0(t), \ph t \ge t_0, \ph u > 0,$

\noindent
where $Q_0(t), \ph P_0(t), \ph F_0(t), \ph L_0(t), \ph M_0(t)$ are  some real-valued continuous functions on $[t_0, +\infty)$. Then every solution $(\phi(t),\psi(t))$ of the system (1.1) with $\phi(t_0) > 0, \ph \psi(t_0) \ge 0$ exists on $[t_0,+\infty)$ and
$$
\phi(t) > 0, \ph \psi(t) \ge 0, \ph t \ge t_0. \eqno (3.33)
$$
}

Proof. Let $(\phi(t), \psi(t))$ be a solution of the system (1.1) with $\phi(t_0) > 0, \ph \psi(t_0) \ge t_0$ and let $[t_0,T)$ be the maximum existence interval for it. We must show that
$$
T = +\infty. \eqno (3.34)
$$
At first show that
$$
\phi(t) > 0, \ph t\in [t_0, T). \eqno (3.35)
$$
Suppose this is not true. Then since $\phi(t_0) > 0$ we have $\phi(t) > 0. \ph t\in [t_0,T_1)$ and
$$
\phi(T_1) = 0 \eqno (3,36)
$$
for some $ T_1 \in (t_0,T)$. By the second and third relation of (2.8) we have
$$
\phi(t) = \phi(t_0)\exp\biggl\{\il{t_0}{t}\Bigl[P(\tau, \phi(\tau),\psi(\tau)) + y(\tau) Q(\tau, \phi(\tau), \psi(\tau))\Bigr]d\tau\biggr\} +
$$
$$
\il{t_0}{t}\exp\biggl\{\il{\tau}{t}\Bigl[P(s, \phi(s),\psi(s)) + y(s) Q(s, \phi(s), \psi(s))\Bigr]d s\biggr\}F(\tau, \phi(\tau), \psi(\tau))d \tau,  \eqno (3.37)
$$
$t\in [t_0,T_1),$
$$
\psi(t) = y(t) \phi(t), \phh t\in [t_0, T_1), \eqno (3.38)
$$
where by by the third relation of (2.8) $y(t)$ is the solution of the Riccati equation
$$
y' + Q(t,\phi(t),\psi(t)) y^2 + \Biggl(B(t,\phi(t),\psi(t)) + \frac{F(t,\phi(t),\psi(t))}{\phi(t)}\Biggr) y - \phantom{aaaaaaaaaaaaaaaaaaaaaaaaaaa}
$$
$$
\phantom{aaaaaaaaaaaaaa}-R(t, \phi(t), \psi(t)) - \frac{G(t,\phi(t),\psi(t))}{\phi(t)} = 0, \phh t\in [t_0,T_1) \eqno (3.39)
$$
with $y(t_0) = \frac{\psi(t_0)}{\phi(t_0)} \ge 0$. Consider the Riccati equation
$$
y' + Q(t,\phi(t),\psi(t)) y^2 + \Biggl(B(t,\phi(t),\psi(t)) + \frac{F(t,\phi(t),\psi(t))}{\phi(t)}\Biggr) y  = 0, \phh t\in [t_0,T_1).  \eqno (3.40)
$$
Obviously $y_1(t)\equiv 0$ is a solution of this equation. Then applying Theorem 2.1 to (3.39) and (3.40) on the basis of the conditions 1) and 5) we conclude that
$$
y(t) \ge 0, \phh t\in [t_0,T_1). \eqno (3.41)
$$
Consider the Riccati equation
$$
y' + Q(t,\phi(t),\psi(t)) y^2 + L_0(t) y - M_0(t) = 0, \phh t\in [t_0, T_1).
$$
Let $y_2(t)$ be a solution of this equation with $y_2(t_0) \ge y(t_0)$. Using Theorem 2.1 to this equation and to Eq. (3.39) and taking into account 1), 4), 5) and (3.41) we conclude that $y_2(t)$ exists on $[t_0, T_1)$ and
$$
y(t) \le y_2(t), \phh t \in [t_0,T_1). \eqno (3.42)
$$
This together with (3.37), (3.39), (3.41), and the conditions 1), 2), 3) implies that \linebreak $\phi(T_1) >~ 0$, which contradicts (3.36). The obtained contradiction proves (3.35). It follows from (3.35) that $y(t)$ is continuable on $[t_0,T)$ as a solution of the equation
$$
y' + Q(t,\phi(t),\psi(t)) y^2 + \Biggl(B(t,\phi(t),\psi(t)) + \frac{F(t,\phi(t),\psi(t))}{\phi(t)}\Biggr) y - \phantom{aaaaaaaaaaaaaaaaaaaaaaaaaaa}
$$
$$
\phantom{aaaaaaaaaaaaaaaaaaaaaaaaaa}-R(t, \phi(t), \psi(t)) - \frac{G(t,\phi(t),\psi(t))}{\phi(t)} = 0, \phh t\in [t_0,T).
$$
and by analogy of the proofs of the inequalities (3.41) and (3.42) it can be proved that
$$
0 \le y(t) \le y_2(t), \phh t\in [t_0,T), \eqno (3.43)
$$
where $y_2(t)$ is a solution of the Riccati equation
$$
y' + Q(t,\phi(t),\psi(t)) y^2 + L_0(t) y - M_0(t) = 0, \phh t\in [t_0, T).
$$
with $y_2(t_0)\ge y(t_0)$. By (2.8) we have

$$
\phi(t) = \phi(t_0)\exp\biggl\{\il{t_0}{t}\Bigl[P(\tau, \phi(\tau),\psi(\tau)) + y(\tau) Q(\tau, \phi(\tau), \psi(\tau))\Bigr]d\tau\biggr\} +
$$
$$
\il{t_0}{t}\exp\biggl\{\il{\tau}{t}\Bigl[P(s, \phi(s),\psi(s)) + y(s) Q(s, \phi(s), \psi(s))\Bigr]d s\biggr\}F(\tau, \phi(\tau), \psi(\tau))d \tau,  \phh t\in [t_0,T).
$$
$$
\psi(t) = y(t) \phi(t), \phh t\in [t_0,T). \eqno (3.44)
$$
This together with the conditions 1), 2), 3) implies that $(\phi(t), \psi(t))$ is bounded on $[t_0,T)$ provided $T < +\infty$. Hence  if $T< +\infty$ then by Lemma 2.2  $[t_0,T)$ cannot be the maximum existence interval for $(\phi(t), \psi(t))$. Therefore (3.34) is valid. From (3.34), (3.35), (3.43) and (3.44) it follows (3.33). The theorem is proved.

{\bf Remark 3.4.} {\it Theorem 3.3 remains valid if we replace its conditions 2), 3)  by the following one
$$
P(t,u,v) + \frac{F(t,u,v)}{u} \le P_0(t), \phh t \ge t_0, \phh u > 0.
$$
which can be substantiated by using in the proof of Theorem 3.3 the first relation of (2.8) instead of the second one of (2.8).}

{\bf Theorem 3.5}. {\it Let $q_0(t), \ph r_0(t), \ph \alpha_0(t), \ph \beta_0(t)$ be real-valued locally integrable on $[a,b] \linebreak (t_0 \le a < b < +\infty)$ functions such that
$$
Q(t,u,v) \ge q_0(t)\ge 0, \ph R(t,u,v) \le r_0(t)\le 0, \ph \alpha_0(t) \le B(t,u,v) \le \beta_0(t), \ph t \in [a,b],
$$
$$
\il{a}{b}\min\biggl[q_0(t)\exp\biggl\{-\il{a}{t}\alpha_0(\tau)d\tau\biggr\}, -r_0(t)\exp\biggl\{\il{a}{t}\beta_0(\tau)d\tau\biggr\}\biggr]d t \ge \pi. \eqno (3.45)
$$
Then the system (1.2) is oscillatory on the interval $[a,b]$.
}

Proof. Let $(\phi(t), \psi(t))$ be a nontrivial solution of the system (1.2) on $[a,b]$. Consider the transcendental system of equations
$$
\sist{\rho\sin\theta = \phi(t) \exp\biggl\{-\il{a}{t}P(\tau,\phi(\tau), \psi(\tau))d \tau\biggr\},}{\rho\cos\theta = \psi(t) \exp\biggl\{-\il{a}{t}S(\tau,\phi(\tau), \psi(\tau))d \tau\biggr\}, \phh t\in [a,b]} \eqno (3.46)
$$
We set
$$
\rho(t) \equiv \sqrt{\phi^2(t)\exp\biggl\{-2\il{a}{t}P(\tau,\phi(\tau),\psi\tau)d\tau\biggr\} + \psi^2(t)\exp\biggl\{-2\il{a}{t}S(\tau,\phi(\tau),\psi\tau)d\tau\biggr\}}, \phantom{aaaaaaaaa}
$$

$$
\theta(t) \equiv \theta_0 + \il{a}{t}\frac{[\phi(\tau) E(\tau)]'\psi(\tau) - \phi(\tau) E(\tau) \psi'(\tau)}{\phi^2(\tau) E^2(\tau) + \psi^2(\tau)} d \tau, \phh t \in [a,b],
$$
where $\theta_0\equiv \sist{\arctan \frac{\phi(a)}{\psi(a)}, \ph \psi(a) \ne 0,}{\pi/2  - \arctan \frac{\psi(a)}{\phi(a)}, \ph \phi(a) \ne 0,} \ph E(t) \equiv \exp\biggl\{\il{a}{t}B(\tau,\phi(\tau),\psi(\tau))d\tau\biggr\}, \linebreak t~ \in [a,b].$
It is not difficult to verify that
$(\rho(t), \theta(t))$ is an absolutely continuously solution of the system (3.46) on $[a,b]$.
Then by (1.2) we have
$$
\rho'(t)\sin \theta(t) + \rho(t) \theta'(t) \cos \theta(t) = Q(t, \phi(t), \psi(t))\exp\biggl\{-\il{a}{t}B(\tau, \phi(\tau),\psi(\tau))d \tau\biggr\}\rho(t) \cos \theta(t),
$$
$$
\rho'(t)\cos \theta(t) - \rho(t) \theta'(t) \sin \theta(t) = R(t, \phi(t), \psi(t))\exp\biggl\{\il{a}{t}B(\tau, \phi(\tau),\psi(\tau))d \tau\biggr\}\rho(t) \sin \theta(t).
$$
almost everywhere on $[a,b]$.  Multiply both sides of these relations by $\cos \theta(t)$ and both sides of the second of the ones by $-\sin \theta(t)$ and summarize. We obtain
$$
\rho(t)\theta'(t) = \biggl[Q(t,\phi(t),\psi(t))\exp\biggl\{-\il{a}{t} B(\tau, \phi(\tau),\psi(\tau))d \tau\biggr\} \cos^2 \theta - \phantom{aaaaaaaaaaaaaaaaaaaaa}
$$
$$
\phantom{aaaaa}- R(t, \phi(t),\psi(t))\exp\biggl\{\il{a}{t}B(\tau, \phi(\tau),\psi(\tau))d \tau\biggr\} \sin^2 \theta\biggr] \rho(t), \phh t \in [a,b].
$$
Since $\rho(t) \ne 0, \ph t \in [a,b]$ from here and from the conditions of the theorem it follows
$$
\theta(b) - \theta(a) \ge \il{a}{b}\min\biggl[q_0(t)\exp\biggl\{-\il{a}{t}\alpha_0(\tau)d\tau\biggr\}, - r_0(t)\exp\biggl\{\il{a}{t}\beta_0(\tau)d \tau\biggr\}\biggr]d t \ge \pi.
$$
This together wit (3.46) implies that the system (1.2) is oscillatory on the interval $[a,b]$.
The theorem is proved.

{\bf Corollary 3.2}. {\it Let the conditions of Theorem 3.5, except of (3.45), be satisfied for \linebreak $b=+\infty$ and let
$$
\il{a}{+\infty}\min\biggl[q_0(t)\exp\biggl\{-\il{a}{t}\alpha_0(\tau)d\tau\biggr\}, -r_0(t)\exp\biggl\{\il{a}{t}\beta_0(\tau)d\tau\biggr\}\biggr]d t = +\infty. \eqno (3.47)
$$
Then the system (1.2) is oscillatory.
}

Proof. If (3.47) is satisfied, then there exists an infinitely large sequence $a_1 < b_1 < a_2 < b_2 < \dots < a_n < b_n < \dots$ such that
$$
\il{a_n}{b_n}\min\biggl[q_0(t)\exp\biggl\{-\il{a_n}{t}\alpha_0(\tau)d\tau\biggr\}, - r_0(t)\exp\biggl\{\il{a_n}{t}\beta_0(\tau)d \tau\biggr\}\biggr]d t \ge \pi, \ph n =1,2,\dots.
$$
Then by Theorem 3.5 the system (1.2) is oscillatory on each of the intervals $[a_n, b_n], \linebreak n = 1,2, \dots$. Hence, the system (1.2) is oscillatory.
The corollary is proved.

Theorem 3.5 as well as Corollary 3.2 are hypothetic (conditional) in the sense that the existence of solutions of the system (1.2) on appropriate intervals $[a,b]$ and $[a,+\infty)$ respectively is only assumed but not proved. In the regular case, when the "coefficients" \linebreak $P, \ph Q, \ph R, \ph S$ have no singularities, this conditionality is omitted   by Theorem 3.1. The next theorem weakens  the hypothetic character of Theorem 3.5 and Corollary 3.2 (to some extent) in the singular case (when the "coefficients"  $P, \ph Q, \ph R, \ph S$ may have  singularities).

{\bf Theorem 3.6}. {\it Let the following conditions be satisfied
$$
P(t,u,v) < 0, \ph Q(t,u,v) \ge 0, \ph R(t,u,v) \le 0, \ph S(t,u,v) < 0
$$
$$
\frac{Q(t,u,v)}{P(t,u,v)} \le M, \phh \frac{R(t,u,v)}{S(t,u,v)} \le \frac{1}{M}
$$
for some $M > 0$. Then for every $\alpha, \beta \in \mathbb{R}$ any solution $(\phi(t), \psi(t))$ of the system (1.2) with $\phi(t_0) = \alpha, \ph \psi(t_0) = \beta$ exists on $[t_0,+\infty)$.
}

Proof. Let $(\phi(t), \psi(t))$ be a solution of the system (1.2) with  $\phi(t_0) = \alpha, \ph \psi(t_0) = \beta$ and let $[t_0,T)$ be the maximum existence interval for $(\phi(t), \psi(t))$. We must show that
$$
T = +\infty.
$$
Assume $\alpha^2 + \beta^2 \ne 0$. Then by (2.1) and (2.2) two cases are possible

\noindent
$I^0$) there exists $t_1 \in [t_0,T)$ such that $\phi(t)\ne 0, \ph \psi(t) \ne 0, \ph t\in [t_1,T)$

\noindent
$II^0$) there exists infinite sequence $t_0 < t_1 < ... < T$ such that
$$
\sist{\phi(t) > 0, \ph t\in [t_1,t_2),}{\psi(t_1)=0, \ph \psi(t) \le 0, \ph t \in [t_1, t_2);}\phantom{aaaaaaaaaaaaaaaaaaaaaaaaaaaaaaaaaaaaaaaaaaaaaa}
$$
$$
\sist{\phi(t_2) = 0, \ph \phi(t) \le 0, \ph t \in [t_2,t_3),}{\psi(t) \le 0, \ph t \in [t_2,t_3);}\phantom{aaaaaaaaaaaaaaaaaaaaaaaa}
$$
$$
\sist{\phi(t)\le 0, \ph t\in [t_3,t_4),}{\psi(t_3)=0, \ph \psi(t) \ge 0, \ph t \in [t_3, t_4);}\phantom{aaaaaaaa}
$$
$$
\phantom{aaaaaaaaa}\sist{\phi(t_4) = 0, \ph \phi(t) \ge 0, \ph t \in [t_4,t_5),}{\psi(t) > 0, \ph t \in [t_4,t_5);}
$$
$$
\phantom{aaaaaaaaaaaaaaaaaaaaaaaaa}\sist{\phi(t) > 0, \ph t\in [t_5,t_6),}{\psi(t_5)=0, \ph \psi(t) \le 0, \ph t \in [t_5, t_6);}
$$
- - - - - - - - - - - -  - - - - - - - -  - - - - - - - - -  - - - - - - - - - - - - - - - - - - - - - - - - -

\noindent
so on (roughly speaking the point $(\phi(t), \psi(t))$ rotates around of the origin of coordinate plane $\phi, \psi$ in the clockwise direction infinitely many times when $t$ varies from $t_1$ to $T$). Assume the case $I^0$)takes place. and $\phi(t_1) . 0, \ph \psi(t_1) > 0$. Then by (2.1) and (2.2) from the conditions of the theorem we obtain respectively
$$
|\phi(t)| = \biggl|\phi(t_1)\exp\biggl\{\il{t_1}{t}P(\tau,\phi(\tau),\psi(\tau))d\tau\biggr\} -
$$
$$
- \il{t_1}{t}\biggl(\exp\biggl\{-\il{t_1}{\tau}P(s,\phi(s),\psi(s))d s\biggr\}\biggr)' \frac{Q(\tau,\phi(\tau),\psi(\tau))}{P(\tau,\phi(\tau),\psi(\tau))}d\tau \exp\biggl\{\il{t_1}{t}P(\tau,\phi(\tau),\psi(\tau))d\tau\biggr\}\biggr| \le
$$
$$
\le  M|\phi(t_1)|\exp\biggl\{\il{t_1}{t}P(\tau,\phi(\tau),\psi(\tau))d\tau\biggr\}\biggl(\exp\biggl\{\il{t_1}{t} -P(\tau,\phi(\tau),\psi(\tau))d\tau\biggr\} -1\biggr) \le M|\phi(t_1)|,
$$
$$
|\psi(t)| \le |\psi(t_1)|, \phh t \in [t_1, T).
$$
Hence $(\phi(t), \psi(t))$ is bounded on $[t_0,T)$ provided $T < +\infty$. By similar way it can be shown the boundedness of  $(\phi(t), \psi(t))$ on  $[t_0,T)$  for $T < +\infty$ in the remaining subcases:

\noindent
$\phi(t_1) > 0, \phh \psi(t_1) < 0;$

\noindent
$\phi(t_1) < 0, \phh \psi(t_1) > 0;$

\noindent
$\phi(t_1) < 0, \phh \psi(t_1) < 0.$

\noindent
Thus by Lemma 2.2  in the case $I^0$) for $T< +\infty$  the interval $[t_0,T)$ cannot be the maximum existence interval for $(\phi(t), \psi(t))$. In the case $II^0$) using relations (2.1)and (2.2) on the basis of the conditions of the theorem we obtain
$$
|\phi(t)| \le |\phi(t_1)|, \phh t\in [t_1,t_2),
$$
$$
|\psi(t)| \le \biggl|\exp\biggl\{\il{t_1}{t}S(\tau, \phi(\tau), \psi(\tau))d \tau\biggr\} \times \phantom{aaaaaaaaaaaaaaaaaaaaaaaaaaaaaaaaaaaaaaaaaaa}
$$
$$
 \times
 \il{t_1}{t}\biggl(\exp\biggl\{-\il{t_1}{\tau}S(s, \phi(s), \psi(s))d s\biggr\}\biggr)'\frac{R(\tau, \phi(\tau), \psi(\tau))}{S(\tau, \phi(\tau), \psi(\tau))} \phi(\tau) d \tau\biggr| \le \frac{|\phi(t_1)|}{M}, \phh t\in [t_1,t_2).
$$
Using these estimates by similar way we can obtain  successively
$$
|\phi(t)| \le |\phi(t_1)|, \phh |\psi(t)| \le \frac{|\phi(t_1)|}{M}, \phh t \in [t_k, t_{k+1}), \phh k =2, 3, .... .
$$
It follows from here that $(\phi(t), \psi(t))$ is bounded on $[t_0,T)$ provided $T < +\infty$. Therefore in the case $II^0$) for $T +\infty$ the interval $[t_0, T)$ cannot be the maximum existence interval of $(\phi(t), \psi(t))$. If $\alpha = \beta = 0$ then either $\phi(t) = \psi(t) \equiv 0$ form a solution of the system, or $\phi^2(t_1) + \psi^2(t_1) \ne 0$ for some $t_1 > t_0$. In the last case by already proven $(\phi(t), \psi(t))$ exists on $[t_1, +\infty)$. Hence   $(\phi(t), \psi(t))$ exists on $[t_0, +\infty)$. The theorem is proved.

{\bf Remark 3.5}. {\it Theorem 3.6 remains valid if we replace in its conditions the relations
$$
\frac{Q(t,u,v)}{P(t,u,v)} \le M, \phh \frac{R(t,u,v)}{S(t,u,v)} \le \frac{1}{M}
$$
by the following ones
$$
\frac{Q(t,u,v)}{S(t,u,v)} \le M, \phh \frac{R(t,u,v)}{P(t,u,v)} \le \frac{1}{M}.
$$
This  can be substantiated    by using in the proof of Theorem 3.6 the relations (2.3) and (2.4) instead of (2.1) and (2.2) respectively.}

\pagebreak

\centerline{\bf References}

\vskip20pt

\noindent
1. D. W. Jordan and P. Smith, Nonlinear ordinary differential equations. Oxford University \linebreak \phantom{a}  Press. Fourth Edition, 2007.

\noindent
2. N. N. Bogolubov, Yu. A. Mitropolski, Asymptotic methods in the theory of nonlinear  \linebreak \phantom{a}  oscillation. Mskow, "Mir'', 1974, pp. 1 - 410.

\noindent
3. Momeni M., Kurakis I., Mosley - Fard M., Shukla P. K. A Van der Pol - Mathieu  \linebreak \phantom{a}  equation for the
dynamics of dust grain charge in dusty plasmas. J. Phys, A: Math.Theor.,  \linebreak \phantom{a}  v. 40, 2007, F473 - F481.

\noindent
4. N. Yamaoka, Oscillation criteria for second-order damped nonlinear differential equations with $p$-Laplacian, J. Math. Annal. Appl. 322 (2007) 932-948.

\noindent
5. M. Kitano and K. Takasi, On a class of second order auasilinear ordinary differential equations, Hiroshima Math. J. 25 (1995), 321-355.

\noindent
6. X. Wang, G. Song, Oscillation Theorems for a Class of Nonlinear Second Order Differential Equations with Damping, Advances in Pure Mathematics, 2013, 3, 226-233.

\noindent
7. 4. G. A. Grigorianm  On two comparison tests for second-order linear ordinary differential \linebreak \phantom{a}
equations, Diff. Urav. 47 (2011), 1225(1240 (in Russian), Diff. Eq. 47 (2011), 1237)
1252 \linebreak \phantom{a} (in English).


\noindent
8. 6. Bellman R. Stability theory of differential equations, Moscow, Izdatelstvo inostrannoj \linebreak \phantom{a} literatury, 1954
(New York, Toronto, London, McGRAW-HILL BOOK COMPANY, \linebreak \phantom{a} INC, 1053).

\noindent
9. 7. Ph. Hartman, Ordinary Differential Equations, SIAM, Classics in Applied Mathematics \linebreak \phantom{a} literature, 1954, 38, Philadelphia 2002.

\vskip20pt

\end{document}